\documentclass[a4paper,10pt]{article}
\usepackage{amsmath}
\usepackage{amssymb}
\usepackage{graphics}
\usepackage{rotating}
\usepackage{cite}
\usepackage{amsfonts}
\usepackage{color}

\newcommand{\field}[1]{\mathbb{#1}}
\newtheorem{definition}{Definition}
\newtheorem{lemma}{Lemma}

\newtheorem{proposition}{Proposition}
\title{Semiweak Cullen-regular quaternionic functions}
\author{Daniel Alay\'{o}n-Solarz}

\begin{document}

\maketitle

\begin{abstract}

We introduce the class of semiweak  Cullen-regular quaternionic functions by interpreting Cullen-regular functions as solutions to an inhomogeneous PDE in terms of the Fueter operator.
\end{abstract}

\section{Introduction} 
This paper is a continuation of \cite{Ala1}, where the class of Cullen-regular functions \cite{GG} is characterized in terms of the Fueter operator. In the Fueter operator framework the Cullen-regular functions look like a solution to a certain inhomogeneous equation. This in turn leads to a integral theorem for the Cullen functions and thus to the concept of weak solution. However this was just barely sketched. Here we adapt some ideas from the theory of Fundamental Solutions and Weak solutions to homogeneous and inhomogeneous PDE  as given in \cite{Tut}.  A corresponding version for this results exists in Clifford Algebras also so they are hardly new. It is interesting to compare the proof hereby given of the Cauchy Representation Theorem for Fueter-regular functions with \cite{Sudbery}. In view of this, most of this paper can be read as a short handbook. Our contribution consists on observing that a version of C-regularity that requires weaker hypothesis than $C^{1}$ can be easily defined.

\section{Solutions to the Fueter operator in the distributional sense}
First we recall that as consequence of the Gauss Theorem in 4 dimensions, one has a Quaternionic Gauss Theorem.

\begin{lemma}
Let $(n_{0},n_{1},n_{2},n_{3})$ be the outward unit normal of  $\partial \Omega$,  any smooth, simple closed hypersurface in the quaternionic space and $\Omega$ its interior, then
\begin{displaymath}
\int_{\partial \Omega}(f_{0}n_{0}+f_{1}n_{1}+f_{2}n_{2}+f_{3}n_{3}) dS_{K} = \int_{\Omega}(\frac{\partial f_{0}}{\partial t}+\frac{\partial f_{1}}{\partial x}+\frac{\partial f_{2}}{\partial y}+\frac{\partial f_{3}}{\partial z}) dV
\end{displaymath}
where $f_{i}$ are differentiable quaternionic functions.
\end{lemma}
and then, a Quaternionic Green Integral Formula, denoting $D_{l},D_{r}$ as the left- and right-Fueter operator:
\begin{lemma}
Let $u$ and $v$ be quaternionic functions defined in a smooth domain $\Omega$ and its border $\partial \Omega$ and, then 
\begin{displaymath}
\int_{\Omega}{(D_{r}u)v + u(D_{l}v)} \ dV = \int_{\partial \Omega} unv \  dS.
\end{displaymath}
\end{lemma}
\textbf{Proof} Let $u$ and $v$ be a pair of differentiable quaternionic functions satisfying the hipothesis, then we set:
\begin{displaymath}
f_{0} := uv,
\end{displaymath}
\begin{displaymath}
f_{1} := uiv,
\end{displaymath}
\begin{displaymath}
f_{2} := ujv,
\end{displaymath}
\begin{displaymath}
f_{3} := ukv,
\end{displaymath}

and apply \textbf{Lemma 1}.

\begin{definition}
Let $f$ is a integrable quaternionic function defined in a bounded domain $\Omega$ with boundary $\partial \Omega$. We say $f$ is  \textbf{weak F-regular} in $\Omega$ if
\begin{displaymath}
\int_{\Omega} (D_{r} \varphi) f  \ dV = 0,
\end{displaymath}
for all test functions $\varphi$ vanishing in $\partial \Omega$.
\end{definition}

\begin{proposition}
Let $f$ be a $C^{1}$ weak F-regular function, then $f$ is F-regular in the classical sense.
\end{proposition}
\textbf{Proof} Suppose $f$ is a $C^{1}$ weak F-regular function. Then we have, setting $u:= \varphi$ a test function and $u:= f$ in \textbf{Lemma 1}
\begin{displaymath}
\int_{\Omega} \varphi (D_{l} f)  \ dV = 0,
\end{displaymath}
for all choices of test functions $\varphi$. In particular if $\varphi$ is real-valued, applying the Fundamental Lemma of Variational Calculus on the components of $D_{l} f$ we conclude that
\begin{displaymath}
D_{l} f = 0.
\end{displaymath}
in $\Omega$.
More generally,
\begin{definition}
Let $f$ and $h$ be integrable quaternionic functions defined in a bounded domain $\Omega$ with boundary $\partial \Omega$. We say $f$ is  \textbf{weak} solution to the inhomogeneous equation 
\begin{displaymath}
D_{l} f = h,
\end{displaymath}
in $\Omega$ if 
\begin{displaymath}
\int_{\Omega} (D_{r} \varphi) f + \varphi h \ dV = 0
\end{displaymath}
for all test functions $\varphi$.
\end{definition}
\begin{proposition}
Let $f$ be a $C^{1}$ weak solution to the inhomogeneous equation, then $f$ is a solution in the classical sense.

\end{proposition}
\textbf{Proof} Let $f$ be a such function. Define a continuous quaternionic function $g$ as
\begin{displaymath}
g:= D_{l}f -h,
\end{displaymath}
then, by hypothesis, for all test functions $\varphi$:
\begin{displaymath}
\int_{\Omega}{(D_{r} \varphi) f   + h \varphi} \ dV  = 0,
\end{displaymath}
now, since for all test functions $\varphi$ we have
\begin{displaymath}
\int_{\Omega}(D_{r} \varphi) f  \ dV= -\int_{\Omega} (D_{l} f) \varphi \  dV,
\end{displaymath}
this implies, using the hipothesis that
\begin{displaymath}
0 = -\int_{\Omega}{(D_{l} f) \varphi \ dV + \int_{\Omega}{  h \varphi}} \ dV =  - \int_{\Omega}  g \varphi \ dV,
\end{displaymath}
which holds in particular if the test function is real. Using the Fundamental Lemma of Variational Calculus on the components of $g$ we see that $g$ must vanish identically in $\Omega$. 
We obtain then that the $C^{1}$ function $f$ satisfy
\begin{displaymath}
D_{l}f = h.
\end{displaymath}

\section{Inverse of the Fueter operator in the distributional sense}
Let $E(q,p)$ denote the following quaternion valued function of the two quaternionic variables $p$ and $q$.
\begin{displaymath}
E(q,p):=  \frac{1}{2 \pi^{2}}\frac{(q-p)^{-1}}{ |q - p|^{2} }.
\end{displaymath}
Note that 
\begin{displaymath}
|E(q,p)| = \frac{1}{2 \pi^2 |q-p|^{3}},
\end{displaymath}
and
\begin{displaymath}
E(q,p) = - E(p,q).
\end{displaymath}
For a fixed quaternion $q$ the function $E$ is analytic for every $p \neq q$ with a isolated weak singularity at $q$. $E$ is left and right F-regular on its domain $\field{H} \backslash \{q\}$:
\begin{displaymath}
D_{l} E = D_{r} E = 0.
\end{displaymath}
\begin{proposition} Let $f$ be a continuous function, $n$ the unit outward normal to the ball with center $p$, then
\begin{displaymath}
\lim_{\epsilon \to 0} \int_{|q-p|= \epsilon} E(q,p)n(q) f(q) dS_{q} = f(p).
\end{displaymath}
\end{proposition}
\textbf{Proof} For every positive $\epsilon$, let $d \mu$ denote the element of surface of the unit 3-shere. Then:
\begin{displaymath}
\lim_{\epsilon \to 0} \int_{|q-p|= \epsilon} E(q,p)n(q) f(q) dS =  \lim_{\epsilon \to 0} \int_{|q-p|= \epsilon}  f(q) d \mu = f(p).
\end{displaymath}
\textbf{Remark} For this limit to hold it suffices that the function $f$ is continuous.
\begin{proposition} Let $\varphi$ be a test function on $\Omega$ and let $p \in \Omega$ be an interior point  then
\begin{displaymath}
\int_{\Omega} {(D_{r} \varphi) E(q,p)} \ dV_{q} = -\varphi(p).
\end{displaymath}
\end{proposition}

\textbf{Proof} Since we want to apply \textbf{Lemma 2} to a function having a isolated singularity at an interior point $p$ of $\Omega$ we need to omit a neighbourhood of $p$. Let $\Omega_{\epsilon} := \Omega \backslash \bar{U_{\epsilon}}$ where $U_{\epsilon}$ denotes the $\epsilon$-neighbourhood of $p$. Then the boundary of $\Omega_{\epsilon}$ consists of $\partial \Omega$ and the $\epsilon$-sphere centered at $p$. Now, applying \textbf{Lemma 2} with $u :=  \varphi$ a test function and $v := E$  we have
\begin{displaymath}
\int_{\Omega_{\epsilon}}{\varphi (D_{l} E(q,p)) + (D_{r} \varphi) E(q,p)} \ dV_{q} =
\end{displaymath}
\begin{displaymath}
 \int_{\partial \Omega} E(q,p) \ n \  \varphi \  dS_{q} +  \int_{|q-p|=\epsilon} E(q,p) \ n(q) \  \varphi \  dS_{q},
\end{displaymath}
and from this we conclude that
\begin{displaymath}
\int_{\Omega_{\epsilon}}{ (D_{r} \varphi) E(q,p)} \ dV_{q} = \int_{|q-p|=\epsilon} E(q,p) \ n(q) \  \varphi(q) \  dS_{q}.
\end{displaymath}
Now, applying the \textbf{Proposition 2}, since $n$ is an \textit{inward} unit normal in the $\epsilon$-sphere:
\begin{displaymath}
\lim_{\epsilon \to 0} \int_{|q-p|=\epsilon} E(q,p) \ n(q) \  \varphi(q) \  dS_{q} = -\varphi(p),
\end{displaymath}
from this we conclude that
\begin{displaymath}
\int_{\Omega} {(D_{r} \varphi) E(q,p)} \ dV_{p} = -\varphi(p).
\end{displaymath}
\begin{proposition}
The function $g$ defined as
\begin{displaymath}
g(p) := -\int_{\Omega} E(q,p) \ h(q)  \ dV_{q},
\end{displaymath} 
is a weak solution to the equation
\begin{displaymath}
D_{l} g = h.
\end{displaymath}
\end{proposition}
\textbf{Proof} Denoting $\Omega$ as domain of the $p$- and $q$-space as $\Omega_{p}$ and $\Omega_{q}$, we have:
\begin{displaymath}
\int_{\Omega_{q}} (D_{r} \varphi) g(q) dV_{q} = \int_{\Omega_{q}} (D_{r} \varphi) \Big (- \int_{\Omega_{p}} E(p,q) \ h(p)  \ dV_{p} \Big) dV_{q} 
\end{displaymath}
Using Fubini's Theorem for weakly singular integrals the right hand side can be written as
\begin{displaymath}
\int_{\Omega_{p}}  \Big ( \int_{\Omega_{q}} (D_{r} \varphi)   E(q,p)  dV_{q} \Big) \ h(p)  \ dV_{p}   = \int_{\Omega_{p}}  -\varphi(p) \ h(p)  \ dV_{p},
\end{displaymath}
and thus
\begin{displaymath}
\int_{\Omega_{q}} (D_{r} \varphi) g(q) + \varphi(q) \ h(q) \ dV_{q} =0. 
\end{displaymath}
\section{Analiticity of Fueter-regular functions}
\begin{proposition}
Let $f$ be a F-regular function  defined in a domain $\Omega$ and its boundary $\partial \Omega$ and $p$ a point in the interior of $\Omega$ then 
\begin{displaymath}
f(p) = \int_{\partial \Omega} E(q,p)n(q) f(q) dS_{q}
\end{displaymath}
and therefore F-regular functions are analytic.
\end{proposition}
\textbf{Proof} Setting $u:= E$ and $v = f$ a F-regular function in \textbf{Lemma 2} we have:
\begin{displaymath}
\int_{\Omega_{\epsilon}}  (D_{r} E(q,p)) f(q) +  E(q,p) (D_{l}(f(q)) \ dV_{q} =
\end{displaymath}
\begin{displaymath}
 \int_{\partial \Omega} E(q,p) n(q) \ f(q) \   \  dS_{q} +  \int_{|q-p|=\epsilon} E(q,p) \ n(q) \  f(q) \  dS_{q}.
\end{displaymath}
Since the left handside vanishes for every $\epsilon$ we have:
\begin{displaymath}
 \int_{\partial \Omega} E(q,p) \ f(q) \   \  dS_{q} +  \int_{|q-p|=\epsilon} E(q,p) \ n(q) \  f(q) \  dS_{q} = 0.
\end{displaymath}
Now taking the limit as $\epsilon$ goes to zero, since the unit normal on the $\epsilon$-sphere is inward, using \textbf{Proposition 2} we have the desired result.

\section{Special case: Semiweak Cullen-regular functions}
Define $\iota$ as
\begin{displaymath}
\iota := \frac{xi+yj+zk}{\sqrt{x^{2}+y^{2}+z^{2}}}.
\end{displaymath}
then we write $p$ as $t + r \iota$, with $r$ the norm of the imaginary part of $p$. \\
Quaternionic functions of one quaternionic variable that are of class $C^{1}$ and null-solutions to the \textit{left-Cullen operator} defined as:
\begin{displaymath}
 \frac{\partial }{\partial t} + \iota \frac{\partial }{\partial r},
\end{displaymath}
will be called \textbf{C-regular}. \\
Note that $\iota$ can be parametrized by spherical coordinates
\begin{displaymath}
\iota = (\cos\alpha \sin\beta, \sin\alpha \sin\beta, \cos \beta).
\end{displaymath}
The left-Fueter operator in this coordinate system has the form:
\begin{displaymath}
D_{l} = \frac{\partial }{\partial t} + \iota \frac{\partial }{\partial r}  - \frac{1}{r} \frac{\partial}{\partial \iota},
\end{displaymath}
for all $p$ that does not lie in the complex subplane of the quaternions given by $t + zk$. Where
\begin{equation*}
\frac{\partial}{\partial \iota} := ({\iota}_{\alpha})^{-1}\frac{\partial}{\partial \alpha} +  ({\iota}_{\beta})^{-1}\frac{\partial}{\partial \beta}.
\end{equation*}
and $\iota_{\alpha}$ and $\iota_{\alpha}$ represent the derivatives of $\iota$ with respect to $\alpha$ and $\beta$ respectively. 

\begin{definition}
Let $f$ be a integrable quaternionic function defined in a bounded domain $\Omega$ with boundary $\partial \Omega$ that admits derivatives in the variables $\alpha$ and $\beta$. We say $f$ is  \textbf{semi-weak C-regular} in $\Omega$ if
\begin{displaymath}
\int_{\Omega}{(D_{r} \varphi) f \ dV = \int_{\Omega}{ (\frac{2 v  \varphi }{r}})} \ dV  
\end{displaymath}
where
\begin{displaymath}
v:= \frac{1}{2} \frac{\partial f}{\partial \iota},
\end{displaymath}
for all test functions $\varphi$ vanishing in $\partial \Omega$.
\end{definition}
\begin{proposition}
Let $f$ be a $C^{1}$ semi-weak C-regular function, then $f$ is C-regular in the classical sense.
\end{proposition}
\textbf{Proof} Applying \textbf{Proposition 2} to $f$ we see that $f$ must satisfy:
\begin{displaymath}
D_{l}f = -\frac{1}{r} \frac{\partial f}{\partial \iota}
\end{displaymath}
now the latter can be written as:
\begin{displaymath}
(\frac{\partial }{\partial t} + \iota \frac{\partial }{\partial r}  - \frac{1}{r} \frac{\partial}{\partial \iota})f = -\frac{1}{r} \frac{\partial f}{\partial \iota}
\end{displaymath}
and thus
\begin{displaymath}
(\frac{\partial }{\partial t} + \iota \frac{\partial }{\partial r})f = 0.
\end{displaymath}
so $f$ is C-regular.
We now show a representation theorem for C-regular functions in the classic sense. Taking $u= E$ and $v = f$ a C-regular function in \textbf{Lemma 2} we have that for in interior point of $p$ of $\Omega$:
\begin{displaymath}
f(p) =  \int_{\Omega}{E(q,p)(\frac{2v}{r})} \ dV + \int_{\partial \Omega} E(q,p)n(q)f(q) \  dS .
\end{displaymath}
\section{Acknowledgments}
I am grateful to C. Vanegas for interesting discussions on this subject. I am also grateful to W. Tutschke for teaching me about the definition, properties and applications of the concept of Fundamental Solutions in PDE's.

%%%%%%%%%%%%%%%%%%%%%%%%%%%%%%%%%%%%%%%%%%%%%%%%%%%%%%%%%%%%%
%                                SECTION V
%%%%%%%%%%%%%%%%%%%%%%%%%%%%%%%%%%%%%%%%%%%%%%%%%%%%%%%%%%%%%%%%%%%%%%%%%%%%%%%

%%%%%%%%%%%%%%%%%%%%%%%%%%%%%%%%%%%%%%%%%%%%%%%%%%%%%%%%%%%%%%%%%%%%%%%%%%%%%%%
%                                SECTION VI
%%%%%%%%%%%%%%%%%%%%%%%%%%%%%%%%%%%%%%%%%%%%%%%%%%%%%%%%%%%%%%%%%%%%%%%%%%%%%%%

%%%%%%%%%%%%%%%%%%%%%%%%%%%%%%%%%%%%%%%%%%%%%%%%%%%%%%%%%%%%%%%%%%%%%%%%%%%%%%%
%                                SECTION VII
%%%%%%%%%%%%%%%%%%%%%%%%%%%%%%%%%%%%%%%%%%%%%%%%%%%%%%%%%%%%%%%%%%%%%%%%%%%%%%%

%%%%%%%%%%%%%%%%%%%%%%%%%%%%%%%%%%%%%%%%%%%%%%%%%%%%%%%%%%%%%%%%%%%%%%%%%%%%%%%
%                                REFERENCES
%%%%%%%%%%%%%%%%%%%%%%%%%%%%%%%%%%%%%%%%%%%%%%%%%%%%%%%%%%%%%%%%%%%%%%%%%%%%%%%

\end{document}